\documentclass[a4paper,11pt,leqno]{amsart}
\usepackage[latin1]{inputenc}
\usepackage{amsmath}
\usepackage{amsfonts}
\usepackage{amssymb}
\usepackage{dsfont}
\usepackage{mathrsfs}
\usepackage{graphicx}
\usepackage{setspace}
\usepackage{fancyhdr}
\usepackage{amsthm}
\usepackage{empheq}
\usepackage{cases}
\usepackage[all]{xy}
\newcommand{\varr}{(M, \metric)}
\newcommand{\soliton}{\pa{M, \metric, X}}
\newcommand{\erre}{\mathbb{R}}

\newcommand{\cinf}{C^{\infty}(M)}

\newcommand{\ricc}{\operatorname{Ricc}}
\newcommand{\diver}{\operatorname{div}}
\newcommand{\hess}{\operatorname{Hess}}
\newcommand{\dist}{\operatorname{dist}}

\newcommand{\ra}{\rightarrow}

\newcommand{\set}[1]{{\left\{#1\right\}}}               
\newcommand{\pa}[1]{{\left(#1\right)}}                  
\newcommand{\sq}[1]{{\left[#1\right]}}                  
\newcommand{\abs}[1]{{\left|#1\right|}}                 
\newcommand{\pair}[1]{\left\langle#1\right\rangle}      

\newcommand{\metric}{\pair{\;,}}                          

\renewcommand{\qedsymbol}{q.e.d.}

\theoremstyle{plain}
\newtheorem{teorema}{Theorem}[section]
\newtheorem{proposiz}[teorema]{Proposition}
\newtheorem{lemma}[teorema]{Lemma}
\newtheorem{coroll}[teorema]{Corollary}
\numberwithin{equation}{section}

\theoremstyle{definition}

\newtheorem*{oss}{Remark}

\theoremstyle{remark}

\title{Articolo}
\author{Paolo Mastrolia}
\address{Dipartimento di Matematica\\
Universit\`a degli Studi di Milano\\
via Saldini 50\\
I-20133 Milano, ITALY} \email{paolo.mastrolia@unimi.it}

\author{Marco Rigoli}
\address{Dipartimento di Matematica\\
Universit\`a degli Studi di Milano\\
via Saldini 50\\
I-20133 Milano, ITALY} \email{marco.rigoli@unimi.it}

\begin{document}
\title{On the Geometry of complete Ricci Solitons}
\maketitle

%
%
%
%

\begin{abstract}
  In this paper we establish three basic equations for a general soliton structure on the Riemannian manifold $\varr$. We then draw some geometric conclusions with the aid of the maximum principle.
\end{abstract}



\begin{section}{Introduction and main results}

Let $\varr$ be an $m$-dimensional, complete, connected Riemannian manifold. A \emph{soliton structure} $\soliton$ on $M$ is the choice (if any) of a smooth vector field $X$ on $M$ and a real constant $\lambda$ such that
\begin{equation}\label{1.1soliton_def}
  \ricc + \frac 12 \,\mathcal{L}_X\metric = \lambda \metric,
\end{equation}
where $\ricc$ denotes the Ricci tensor of the metric $\metric$ on $M$ and $\mathcal{L}_X\metric$ is the Lie derivative of this latter in the direction of $X$. In what follows we shall refer to $\lambda$ as to the \emph{soliton constant}. The soliton is called \emph{expanding}, \emph{steady} or \emph{shrinking} if, respectively, $\lambda < 0$, $\lambda = 0$ or $\lambda >0$. If $X$ is the gradient of a potential $f \in \cinf$, then \eqref{1.1soliton_def} takes the form
\begin{equation}\label{1.2gradientsoliton_def}
  \ricc + \operatorname{Hess}(f) = \lambda \metric,
\end{equation}
and the Ricci soliton is called a \emph{gradient Ricci soliton}. Both equations \eqref{1.1soliton_def} and \eqref{1.2gradientsoliton_def} can be considered as perturbations of the Einstein equation
\[
\ricc = \lambda\metric
\]
and reduce to this latter in case $X$ or $\nabla f$ are Killing vector fields. When $X=0$ or $f$ is constant we call the underlying Einstein manifold a \emph{trivial} Ricci soliton.

Since the appearance of the seminal works of R. Hamilton, \cite{hamilton}, and G. Perelman, \cite{perelman}, the study of gradient Ricci solitons has become the subject of a rapidly increasing investigation directed mainly towards two goals, \emph{classification} and \emph{triviality}; among the enormous literature on the subject we only quote, as a few examples, the papers \cite{petwylie1}, \cite{petwylie2}, \cite{petwylie3}, \cite{prims}, \cite{prrims}, \cite{ELnM}.

In this paper we focus our attention on the more general case of \emph{Ricci solitons}, that is, when $X$ is not necessarily the gradient of a potential $f$. A first important difference is that, in the present case, we cannot make use of the weighted manifold structure $\pa{M, \metric, e^{-f}\textrm{d}vol}$ which naturally arises when dealing with gradient solitons. The same applies for related concepts such as the Bakry-Emery Ricci tensor, giving rise to weighted volume estimates, or the weak maximum principle for the diffusion operator $\Delta_f$ (the ``$f$-Laplacian''), acting on $u \in C^2(M)$ by
\[
\Delta_f u = \Delta u - \pair{\nabla f, \nabla u},
\]
that we considered in a previous investigation (see \cite{prrims} for details). Thus our assumptions and techniques have to rest on the original Riemannian structure of $M$. This restricts the applicability of the technical tools we used in \cite{prrims}; nevertheless, we are still able to obtain some stringent geometric conclusions as those we are going to describe in a shortwhile.

From now on we fix an origin $o \in M$ and let $r(x) = \dist(x, o)$. We set $B_r$ and $\partial B_r$ to denote respectively the geodesic ball of radius $r$ centered at $o$ and its boundary.

To state our first result, we recall that, given a Schrödinger type operator $L = \Delta+q(x)$ for some $q(x) \in C^0(M)$, the \emph{spectral radius} of $L$ on $M$ is defined via the Rayleigh characterization by
\[
\lambda_1^L(M) =  \inf_{\begin{array}{cc}\scriptstyle{\varphi \in C^{\infty}_0}(M) \\ \scriptstyle{\varphi \not \equiv 0}\end{array}} \frac{\int_M{\abs{\nabla\varphi}^2 -q(x)\varphi^2}}{\int_M \varphi^2}.
\]

\begin{teorema}\label{TheoremANonextLp}
  Let $\varr$ be a complete manifold with Ricci tensor satisfying
  \begin{equation}\label{1.19CondRiccileq1.3}
  \ricc \leq \frac 12 a(x) \metric
  \end{equation}
  for some $a(x) \in C^{0, \alpha}(M)$ with $0 \leq \alpha <1$.
  Assume that, for some $H \geq 1$,
  \begin{equation}\label{1.4newFirstEigenLgeq0}
  \lambda^{L_H}_1(M) \geq 0
  \end{equation}
   where $L_H = \Delta +Ha(x)$. If there exists a soliton structures $\soliton$ on $\varr$ with $X \not \equiv 0$ satisfying
  \begin{equation}\label{1.5newIntConditionX}
     \set{\int_{\partial B_r}\abs{X}^{4(\beta+1)}}^{-1} \not \in L^1(+\infty)
  \end{equation}
  for some $0 \leq \beta \leq H-1$, then $X$ is a parallel field and $\varr$ is Einstein. Furthermore, the simply connected universal cover of $M$ is a warped product $\pa{\erre \times_c P, h}$ with $c = \abs{X}$, $h = \textrm{d}t^2 + c g$ and $\pa{P, \,g}$ is Einstein.
\end{teorema}
\begin{oss}
  If $\abs{X} \in L^{4(\beta +1)}(M)$ then \eqref{1.5newIntConditionX} is satisfied. Furthermore, if $X \equiv 0$ then $\varr$ is trivially Einstein.
\end{oss}

The following is a ``second version'' of Theorem \ref{TheoremANonextLp}.
\begin{teorema}\label{TheoremA'SobolevNonextLp}
 Let $\varr$ be a complete manifold and assume, for some $0\leq \alpha <1$, the validity on $M$ of the Sobolev-Poincaré  inequality
 \begin{equation}\label{1.5'PoincSobolev}
   \int_M \abs{\nabla \varphi}^2 \geq S(\alpha)^{-1} \set{\abs{\varphi}^{\frac{2}{1-\alpha}}}^{1-\alpha}
 \end{equation}
 for each $\varphi \in C^{\infty}_0(M)$ with a positive constant $S(\alpha)$. Suppose
  \begin{equation*}\tag{\ref{1.19CondRiccileq1.3}}
  \ricc \leq \frac 12 a(x) \metric
  \end{equation*}
  for some $a(x) \in C^{0, \alpha}(M)$ with $0 \leq \alpha <1$; let $\sigma > \frac 12$ and assume that
  \begin{equation}\label{1.5''condaplus}
    \|a_+(x)\|_{L^{\frac 1\alpha}(M)} < 4 \frac{2\sigma-1}{\sigma^2}\frac{1}{S(\alpha)}.
  \end{equation}
  Then there are no Ricci solitons $\soliton$ on $\varr$ satisfying $X \not \equiv 0$ and
  \begin{equation}\label{1.5'''condIntX}
    \int_{B_r} \abs{X}^{2\sigma} = o\pa{r^2} \quad \text{as}\,\, r \ra +\infty.
  \end{equation}
\end{teorema}

 We then consider a general triviality result:

\begin{teorema}\label{TH_triv_expander}
  Let $\varr$ be a complete manifold with Ricci tensor satisfying
  \begin{equation}\label{Ricc_leq_a}
    \ricc \leq -\frac{B}{\pa{1+r(x)}^\mu}
  \end{equation}
   for some constant $B>0$ and let $\sigma, \mu \in \erre$ satisfy
  \begin{equation}
    \sigma \geq 0, \qquad \sigma+\mu < 2.
  \end{equation}
  Suppose that
  \begin{equation}\label{liminfLogVol}
    \liminf_{r(x)\ra +\infty} \frac{\log \operatorname{vol}(B_r)}{r(x)^{2-\sigma-\mu}} =d_0 <+\infty.
  \end{equation}
   Then there are no  Ricci solitons $\soliton$  on $M$ satisfying $X \not \equiv 0$ and
  \begin{equation}\label{limsupAbsX}
    \limsup_{r(x)\ra +\infty} \frac{\abs{X}^2}{r(x)^\sigma} \begin{cases}
      = 0, &0<\sigma<2, \\ <+\infty, &\sigma=0.
    \end{cases}
  \end{equation}

\end{teorema}

Our next result relates the soliton constant $\lambda$ with the infimum $S_*$ of the scalar curvature of the manifold $\varr$. We recall that this latter is said to \emph{satisfy the Omori-Yau maximum principle} if for each $u \in C^2(M)$ with $u^* = \sup_M u < +\infty$ there exists a sequence $\set{x_k} \subset M$ such that
\[
(i) \, u(x_k) > u^* - \frac 1k, \quad (ii) \, \abs{\nabla u(x_k)} < \frac 1k, \quad (iii) \, \Delta u(x_k) < \frac 1k
\]
for each $k \in \mathds{N}$. Conditions to insure the validity of the Omori-Yau maximum principle are discussed in \cite{prsmax}.
\begin{teorema}\label{TeorBScalarEstimates}
  Let $\varr$ be a complete manifold of dimension $m$ with scalar curvature $S(x)$ and satisfying the Omori-Yau maximum principle. Let $\soliton$ be a Ricci soliton on $\varr$ with soliton constant $\lambda$. Assume
  \begin{equation}\label{eq1.6supX}
    \abs{X}^* = \sup_M \abs{X} < +\infty.
  \end{equation}
 Let
  \[
  S_* = \inf_M S.
  \]
  \begin{itemize}
    \item [(i)] \, If $\lambda <0$ then $m\lambda \leq S_* \leq 0$. Furthermore, if $S(x_0) = S_*=m\lambda$ for some $x_0 \in M$, then $\varr$ is Einstein and $X$ is a Killing field, while if $S(x_0) = S_*=0$ for some $x_0 \in M$, then $\varr$ is Ricci flat and $X$ is a homothetic vector field.
    \item [(ii)] If $\lambda =0$ then $S_*= 0$. Furthermore, if $S(x_0) = S_*=0$ for some $x_0 \in M$, then $\varr$ is Ricci flat and $X$ is a Killing field.
    \item [(iii)] If $\lambda >0$ then $0 \leq S_* \leq m\lambda$. Furthermore, if $S(x_0) = S_*=0$ for some $x_0 \in M$, then $\varr$ is Ricci flat and $X$ is a homothetic vector field, while if $S(x_0) = S_*=m\lambda$ for some $x_0 \in M$, then $\varr$ is compact, Einstein and $X$ is a Killing field.
  \end{itemize}
\end{teorema}
\begin{oss}
  In case $X = \nabla f$, that is, the soliton is a gradient Ricci soliton and $\varr$ is Einstein, a complete classification is given in Theorem 1.3 of \cite{prrims}.
\end{oss}
From Theorem \ref{TeorBScalarEstimates} we immediately obtain
\begin{coroll}\label{CO_1.4}
  Let $\varr$ be a complete manifold with scalar curvature $S$ such that
  \[
  S_* = \inf_M S(x) < 0 \qquad (\text{resp. } >0)
  \]
  and Ricci tensor satisfying
  \begin{equation}\label{1.10_Ricci_below}
  \ricc \geq -(m-1)B^2\pa{1 + r^2}
  \end{equation}
  for some $B \geq 0$. Then $\varr$ does not support any shrinking or steady (resp., expanding or steady) soliton with
  \[
  \abs{X}^* = \sup_M \abs{X} < +\infty.
  \]
\end{coroll}

  Note that \eqref{1.10_Ricci_below} implies the validity of the Omori-Yau maximum principle on any complete manifold $\varr$. However, the validity of this latter is guaranteed also in other circumstances such as, for instance, those of the next

\begin{coroll}\label{CO_1.4'}
  Let $\varr$ be a complete manifold admitting a shrinking or steady soliton structure $\soliton$ with $\abs{X}^* = \sup_M \abs{X} < +\infty$. Then any proper minimal immersion of $M$ into $\erre^n$, with $n>m=\dim M$, is totally geodesic.
  \end{coroll}

In our last result we consider the conformally flat case.

\begin{teorema}\label{TH_1.5_Okumura}
  Let $\varr$ be a complete manifold of dimension $m \geq 3$ with scalar curvature $S(x)$, trace-free Ricci tensor $T$ and satisfying the Omori-Yau maximum principle. Assume that $\varr$ is conformally flat and that
  \[
  S^* = \sup_M S(x) < +\infty.
  \]
  Let $\soliton$ be a Ricci soliton on $\varr$ with soliton constant $\lambda$ and $\abs{X}^* = \sup_M \abs{X} < +\infty$. Then either $\varr$ is of constant sectional curvature or $\abs{T}^* = \sup_M \abs{T}$ satisfies
  \begin{equation}\label{1.11_Okumura_est}
    \abs{T}^* \geq \frac 12 \pa{\sqrt{m(m-1)} \lambda - S^* \frac{m-2}{\sqrt{m(m-1)}}}.
  \end{equation}
\end{teorema}

\end{section}

\begin{section}{Preliminary results}
The proof of our results rests on three interesting formulas. The first (eq. \eqref{2.9BochnerNGSol}) is due to Bochner (at least in case $X$ is a Killing field), as reported in \cite{petersen}, page 191. The remaining two, that is equations \eqref{2.25DeltaScalar} and \eqref{2.53_Delta_absTsquared}, have been found for gradient Ricci solitons in \cite{prrims} but are new and in fact unexpected, at least for us, in the present generality. In what follows, to perform computations, we shall use the method of the moving frame referring to a local orthonormal coframe $\set{\theta^i}$ for the metric and corresponding Levi-Civita connection and curvature forms, indicated respectively with $\set{\theta^i_j}$ and $\set{\Theta^i_j}$, $1 \leq i, j, \ldots \leq m = \dim M$. The Einstein summation convention will be in force throughout.

The following generalized version of the Bochner formula is probably well known; we include a proof here for the sake of completeness.
\begin{lemma}\label{LE2.1generalIdentity}(Generalized Bochner formula)
  Let $Y$ be a vector field on $M$. Then
   \begin{equation}\label{2.2generalidentityY}
      \diver\pa{\mathcal{L}_Y\metric}(Y) = \frac 12 \Delta \abs{Y}^2 - \abs{\nabla Y}^2 +\ricc\pa{Y, Y} + \nabla_Y \pa{\diver Y},
    \end{equation}
    where $\mathcal{L}_Y\metric$ is the Lie derivative of the metric in the direction of $Y$.
    \begin{proof}
      Let $\set{e_i}$ be the o.n. frame dual to $\set{\theta^i}$. Then
      \[
      Y = Y^ie_i = Y_ie_i
      \]
      and setting $Y_{ij}$ for the coefficients of the covariant derivative $\nabla Y$ of $Y$ we have
      \begin{equation}\label{2.3covDerivY}
        Y_{ij}\theta^j = \textrm{d}Y_i - Y_k\theta^k_i.
      \end{equation}
      Differentiating \eqref{2.3covDerivY}, using the definition of covariant derivative, the structure equations
      \[
      \textrm{d}\theta^i = -\theta^i_j \wedge \theta^j, \qquad \textrm{d}\theta^i_j = -\theta^i_k \wedge \theta^k_j + \Theta^i_j
      \]
      and the components $R^i_{jkt}$ of the Riemann curvature tensor defined by
      \[
      \Theta^i_j = \frac 12 R^i_{jkt} \theta^k \wedge \theta^t
      \]
      we obtain
      \[
      Y_{ikj} \theta^j \wedge \theta^k = -\frac 12Y_t R^t_{ijk} \theta^j \wedge \theta^k.
      \]
      Thus, inverting the indexes $k$ and $j$
      \[
      Y_{ijk} \theta^k \wedge \theta^j = -\frac 12Y_t R^t_{ikj} \theta^k \wedge \theta^j.
      \]
      Comparing these last two equations we deduce
      \begin{equation}\label{2.4commutationY2nd3rd}
        Y_{ijk} - Y_{ikj} = Y_t R^t_{ijk}.
      \end{equation}
      Since
      \[
      \mathcal{L}_Y\metric = \pa{Y_{ik} + Y_{ki}} \theta^i \otimes \theta^k
      \]
      we have
      \begin{equation}\label{2.5diver}
         \diver\pa{\mathcal{L}_Y\metric}(Y)= Y_iY_{ikk} + Y_iY_{kik}.
      \end{equation}
      From the commutation relation \eqref{2.4commutationY2nd3rd} tracing with respect to $i$ and $k$ we obtain
      \begin{equation}\label{2.4'}
      Y_{kik} = Y_{kki} + Y_tR^t_{kik} = Y_{kki} + Y_tR_{ti},
      \end{equation}
      where, as usual, with $R_{ti}$ we have indicated the components of the Ricci tensor. Thus
      \begin{equation}\label{2.6diverplusRicc}
        Y_i Y_{kik} = \nabla_Y \pa{\diver Y} + \ricc\pa{Y, Y}.
      \end{equation}
      On the other hand, from $\abs{Y}^2 = Y_iY_i$ we deduce
      \[
      \textrm{d}\abs{Y}^2 = 2Y_iY_{ik} \theta^k
      \]
      and
      \[
      \Delta\abs{Y}^2 = 2Y_{ik}Y_{ik} + 2 Y_iY_{ikk},
      \]
      or, in other words,
      \begin{equation}\label{2.7DeltaabsYsq}
        \frac 12 \Delta\abs{Y}^2 = \abs{\nabla Y}^2 + Y_iY_{ikk}.
      \end{equation}
      Substituting \eqref{2.6diverplusRicc} and \eqref{2.7DeltaabsYsq} into \eqref{2.5diver} we immediately obtain \eqref{2.2generalidentityY}.
    \end{proof}
\end{lemma}

\begin{oss}
  If $Y=\nabla f$, $f \in \cinf$, then \eqref{2.2generalidentityY} can be rewritten as
  \begin{align*}
     \frac 12 \Delta \abs{\nabla f}^2 &= \abs{\hess(f)}^2 -\ricc\pa{\nabla f, \nabla f} - \pair{\nabla \Delta f, \nabla f} \\ &+\diver\pa{2 \hess(f)}\pa{\nabla f} \\ &=  \abs{\hess(f)}^2 +\ricc\pa{\nabla f, \nabla f} + \pair{\nabla \Delta f, \nabla f},
  \end{align*}
  that is, \eqref{2.2generalidentityY} in this case is the classical Bochner's formula.
\end{oss}

\begin{proposiz}\label{PR2.8MainFormula}
   Let $\soliton$ be a Ricci soliton on $\varr$. Then
  \begin{equation}\label{2.9BochnerNGSol}
    \frac 12 \Delta \abs{X}^2 =  \abs{\nabla X}^2 -\ricc\pa{X, X}.
  \end{equation}
  \begin{proof}
    We trace the soliton equation \eqref{1.1soliton_def} to obtain
    \begin{equation*}
      S + \diver X = m\lambda,
    \end{equation*}
     and from here we deduce
     \begin{equation}\label{2.10nablaS}
       \nabla S = - \nabla \diver X.
     \end{equation}
     On the other hand, contracting twice the second Bianchi's identities we have the well-known formula
     \begin{equation}\label{2.11nablaSanddiverRicc}
       \nabla S = 2\diver \ricc.
     \end{equation}
     Thus, comparing \eqref{2.10nablaS} and \eqref{2.11nablaSanddiverRicc},
     \begin{equation}\label{2.12nabladiverX}
        \nabla \diver X = -2 \diver \ricc.
     \end{equation}
     Now taking the divergence of \eqref{1.1soliton_def} and using the fact that $\diver\pa{\lambda \metric}=0$ we obtain
     \[
     \diver\pa{\mathcal{L}_X\metric} = -2 \diver \ricc
     \]
     and \eqref{2.12nabladiverX} yields
     \[
     \nabla \diver X = \diver\pa{\mathcal{L}_X\metric}.
     \]
     In particular,
     \begin{equation}\label{2.13}
       \nabla_X\diver X = \diver\pa{\mathcal{L}_X\metric}\pa{X}.
     \end{equation}
     Thus applying \eqref{2.2generalidentityY} of Lemma \ref{LE2.1generalIdentity} we immediately get \eqref{2.9BochnerNGSol}.
     \end{proof}
\end{proposiz}
\begin{oss}
  In case $X=\nabla f$, that is, the soliton is a gradient soliton, \eqref{2.9BochnerNGSol} becomes
  \[
   \frac 12 \Delta \abs{\nabla f}^2 = \abs{\hess(f)}^2 -\ricc \pa{\nabla f, \nabla f}.
  \]
  Then, using the identity
  \[
  \frac 12 \pair{\nabla f, \nabla \abs{\nabla f}^2} = \hess\pa{f}\pa{\nabla f, \nabla f}
  \]
  and the gradient Ricci soliton equation \eqref{1.2gradientsoliton_def}, we deduce
  \[
   \frac 12 \Delta \abs{\nabla f}^2 = \frac 12 \pair{\nabla f, \nabla \abs{\nabla f}^2} + \abs{\hess(f)}^2 -\lambda\abs{\nabla f}^2.
  \]
  This latter, with the aid of the diffusion operator $\Delta_f$ of the Introduction, can be written as
  \begin{equation}\label{2.14DeltafGradSol}
    \frac 12 \Delta_f \abs{\nabla f}^2 = \abs{\hess(f)}^2 -\lambda\abs{\nabla f}^2.
  \end{equation}
  Formula \eqref{2.14DeltafGradSol} has often appeared in the recent literature on gradient Ricci solitons: see for instance \cite{ELnM}, \cite{prrims} and the references therein.
\end{oss}

Before proceeding to the next proposition we need to determine some further ``commutation relations''.

\begin{lemma}\label{LE2.15Commut4thYtraced}
  Let $Y$ be a vector field on $M$. Then
  \begin{equation}\label{2.16Commut4thYtraced}
    Y_{tkkt} - Y_{kktt} = \frac 12 \pair{\nabla S, Y} + \frac 12 \operatorname{tr}\pa{\mathcal{L}_Y\metric \circ \ricc}.
  \end{equation}
  \begin{proof}
    We start from the commutation relations \eqref{2.4commutationY2nd3rd}. By taking covariant derivative we deduce
    \begin{equation}\label{2.17covDerivY3rdComm}
      Y_{ijkt} - Y_{ikjt} = Y_{st}R_{sijk} + Y_sR_{sijk, \,t}.
    \end{equation}
    Next we recall that, by definition of covariant derivative,
    \begin{equation}\label{2.18covDerivY3rd}
      Y_{ijk} \theta^k = \textrm{d}Y_{ij} - Y_{tj} \theta^t_i - Y_{it} \theta^t_j.
    \end{equation}
    Thus differentiating both members of \eqref{2.18covDerivY3rd}, using the structure equations and \eqref{2.18covDerivY3rd} itself, we arrive at
    \[
    Y_{ijkl} \theta^l \wedge \theta^k = -\frac 12 \pa{Y_{tj}R_{tilk} + Y_{it}R_{tjlk}} \theta^l \wedge \theta^k,
    \]
    from which, inverting $k$ with $l$ and summing up, we deduce
    \begin{equation}\label{2.19}
      Y_{ijkl} - Y_{ijlk} = Y_{tj}R_{tikl} + Y_{it}R_{tjkl}.
    \end{equation}
    Now, \eqref{2.16Commut4thYtraced} follows immediately from \eqref{2.17covDerivY3rdComm}, \eqref{2.19}, \eqref{2.11nablaSanddiverRicc} and tracing.
      \end{proof}
\end{lemma}

 For later use we also recall:
 \begin{lemma}
   For the Ricci tensor we have
   \begin{equation}\label{2.21Ricci12}
     R_{ij, \,k} = R_{ji, \,k};
   \end{equation}
   \begin{equation}\label{2.22Ricci23}
     R_{ij, \,k} - R_{ik, \,j} = - R_{tijk, \,t};
   \end{equation}
   \begin{equation}\label{2.23Ricci34}
     R_{ij, \,kl} - R_{ij, \,lk} = R_{is}R_{sjkl} + R_{js}R_{sikl}.
   \end{equation}
   \begin{proof}
     \eqref{2.21Ricci12} is obvious. \eqref{2.22Ricci23} follows from the second Bianchi's identities, while \eqref{2.23Ricci34} can be obtained with the same methods used in the proofs of Lemmas \ref{LE2.1generalIdentity} and \ref{LE2.15Commut4thYtraced}.
   \end{proof}
 \end{lemma}


We are now ready to prove
\begin{proposiz}\label{PR2.24_DeltaScalar}
  Let $\soliton$ be a Ricci soliton with soliton constant $\lambda$ on $\varr$ and let $S(x)$ be the scalar curvature. Then
  \begin{align}\label{2.25DeltaScalar}
    \frac 12 \Delta S &= \frac 12 \pair{\nabla S, X} + \lambda S - \abs{\ricc}^2  \\ \nonumber &= \frac 12 \pair{\nabla S, X} + \lambda S - \frac{S^2}{m} - \abs{\ricc - \frac Sm\metric}^2.
  \end{align}
  \begin{proof}
    We start from the soliton equation \eqref{1.1soliton_def}, which in components reads
    \begin{equation}\label{(1)componentsSolEq}
      R_{ij} + \frac 12 \pa{X_{ij} + X_{ji}} = \lambda \delta_{ij}.
    \end{equation}
    Differentiating \eqref{(1)componentsSolEq} we get
     \begin{equation}\label{(2)DiffcomponentsSolEq}
      R_{ij, \,k} = - \frac 12 \pa{X_{ijk} + X_{jik}}.
    \end{equation}
    From \eqref{2.4commutationY2nd3rd} applied to $X$ and \eqref{(2)DiffcomponentsSolEq} we obtain
    \begin{equation}\label{(4)commutationRiccX}
      2\pa{R_{ik, \,j}-R_{jk, \,i}} = X_{jki}-X_{ikj}+X_tR_{tkji}.
    \end{equation}
    Taking covariant derivatives we deduce the further commutation relation
    \begin{equation}\label{(5)CommutRelCovDerRicc}
        2\pa{R_{ik, \,jt}-R_{jk, \,it}} = X_{jkit}-X_{ikjt}+X_{st}R_{skji}+X_{s}R_{skji, \,t}.
    \end{equation}
    Contracting \eqref{2.23Ricci34} with respect to $i$ and $l$ we get
    \begin{equation}\label{(9)CommutRelCovDerRicc34}
       R_{lj, \,kl} - R_{lj, \,lk} = R_{ls}R_{sjkl} + R_{js}R_{sk};
    \end{equation}
    note that, from \eqref{2.22Ricci23}, we deduce
    \begin{equation}\label{(10)XsRicc}
      X_sR_{ijks, \,j} = X_sR_{skji, \,j} = X_sR_{ik, \,s} - X_sR_{is, \,k}.
    \end{equation}
    Now, from \eqref{(5)CommutRelCovDerRicc},
    \begin{equation}\label{(11)DeltaRicc}
      \Delta R_{ik} = R_{ik, \,jj} = R_{jk, \,ij} + \frac 12 \pa{X_{jkij}-X_{ikjj}+X_{sj}R_{skji}+X_{s}R_{skji, \,j}}.
    \end{equation}
    Inserting \eqref{(10)XsRicc} into \eqref{(11)DeltaRicc} we get
    \begin{equation}\label{(13)DeltaRiccRicc}
      \Delta R_{ik} =  R_{jk, \,ij} + \frac 12 \pa{X_{jkij}-X_{ikjj}} + \frac 12 X_{sj}R_{skji}+  \frac 12 X_sR_{ik, \,s} - \frac 12 X_sR_{is, \,k} .
    \end{equation}
    Next, we rewrite \eqref{2.23Ricci34} in the form
    \begin{equation}\label{(14)Ricci34Modif}
       R_{jk, \,ij} = R_{jk, \,ji} + R_{js}R_{skij} + R_{ks}R_{si}
    \end{equation}
    and we insert \eqref{(14)Ricci34Modif} into \eqref{(13)DeltaRiccRicc} to get
    \begin{align}\label{(15)}
       \Delta R_{ik} &=  \frac 12 X_{sj}R_{skji} + \frac 12 X_sR_{ik, \,s} - \frac 12 X_sR_{is, \,k} \\ \nonumber &+  R_{jk, \,ji} + R_{js}R_{skij} + R_{ks}R_{si} + \frac 12 \pa{X_{jkij}-X_{ikjj}}.
    \end{align}
    From the second Bianchi's identities we recall that
    \begin{equation}\label{(16)}
      2 R_{ik, \,i} = S_k,
    \end{equation}
    so that
    \begin{equation}\label{(17)}
      R_{jk, \,ji} = \frac 12 S_{ki}.
    \end{equation}
    Using \eqref{(17)} into \eqref{(15)} yields
     \begin{align}\label{(18)}
       \Delta R_{ik} &= \frac12 S_{ki} + \frac 12 X_{sj}R_{skji} + \frac 12 X_sR_{ik, \,s} - \frac 12 X_sR_{is, \,k} \\ \nonumber &+ R_{js}R_{skij} + R_{ks}R_{si} + \frac 12 \pa{X_{jkij}-X_{ikjj}}.
    \end{align}
    Note that from the soliton equation \eqref{1.1soliton_def}
    \begin{equation}\label{(19)ConseqSolEq}
      X_{sj} = -X_{js} + 2\lambda \delta_{sj} - 2R_{sj},
    \end{equation}
    and therefore
    \begin{equation}\label{(20)ConseqSolEq2}
      \frac 12X_{sj}R_{skji} = -\frac 12 X_{js}R_{skji} + \lambda R_{ki} - R_{sj}R_{skji}.
    \end{equation}
    Substituting \eqref{(20)ConseqSolEq2} into \eqref{(18)} gives
     \begin{align}\label{(21)}
       \Delta R_{ik} &= \frac12 S_{ki} -\frac 12 X_{js}R_{skji} + \lambda R_{ki} - 2R_{sj}R_{skji} + \\ \nonumber &+\frac 12 X_sR_{ik, \,s} - \frac 12 X_sR_{is, \,k}  + R_{ks}R_{si} + \\ \nonumber &+\frac 12 \pa{X_{jkij}-X_{ikjj}}.
    \end{align}
    We trace \eqref{(21)} with respect to $i$ and $k$ and use the relation
    \[
    \frac 12 X_tR_{kk, \,t} - \frac 12 X_tR_{kt, \,k} = \frac 14 X_tR_{kk, \,t}= \frac 14 X_tS_t
    \]
    so that
     \begin{align}\label{(22)DeltaS}
       \frac 12\Delta S &= \lambda S  -\frac 12 X_{js}R_{sj}  - R_{sj}R_{sj}+ \\ \nonumber &+\frac 12 X_tR_{kk, \,t} - \frac 12 X_tR_{kt, \,k}  +\frac 12 \pa{X_{jkkj}-X_{kkjj}} \\ \nonumber &= \lambda S - \abs{\ricc}^2 -\frac 12 X_{js}R_{sj} + \frac 12 \pa{X_{tkkt}-X_{kktt}} + \frac 14 X_tS_t.
    \end{align}
    Now we apply Lemma \ref{LE2.15Commut4thYtraced} to $X$ and from \eqref{(22)DeltaS} we immediately obtain the desired result.
  \end{proof}
\end{proposiz}
\begin{oss}
  In case $X=\nabla f$, $f \in \cinf$, that is, the soliton is a gradient soliton, we can rewrite \eqref{2.25DeltaScalar} in the form
  \[
  \frac 12 \Delta_f S = \lambda S - \abs{\ricc}^2,
  \]
  which is formula (2.15) (with $\lambda$ constant) of Lemma 2.3 in \cite{prrims}.
\end{oss}

Our aim is now to compute $\Delta\abs{T}^2$, where $T$ is the  traceless Ricci  tensor, that is,
\begin{equation}\label{2.38_traceless_Ricci}
  T_{ij} = R_{ij} - \frac Sm \delta_{ij}.
\end{equation}
Thus
\begin{equation}\label{2.39_traceless_Ricci_norm}
  \abs{T}^2 = \abs{\ricc}^2 - \frac{S^2}{m},
\end{equation}
and it follows that
\begin{align}\label{2.40_Delta_abssquared_T}
  \Delta \abs{T}^2 &= \Delta\abs{\ricc}^2-\frac 1m \Delta S^2= \\ \nonumber &= 2\abs{\nabla \ricc}^2 + 2 \pair{\ricc, \Delta \ricc} - \frac 2m S\Delta S - \frac 2m \abs{\nabla S}^2.
\end{align}
We have
\begin{proposiz}\label{PR_2.41_Delta_abssquared_T}
  Let $\soliton$ be a Ricci soliton with soliton constant $\lambda$ on $\varr$ and let $S(x)$ be the scalar curvature. Then,
  \begin{align}\label{2.42}
    \frac 12 \Delta \abs{T}^2  &= \frac 12 \pair{\nabla\abs{T}^2, X } + \abs{\nabla T}^2 + 2\lambda \abs{T}^2 + \\ \nonumber &+ \frac 2m \pa{\frac{S^2}{m} + \abs{T}^2}S + 2 R_{ik}R_{sj}R_{skij}.
  \end{align}
  \begin{proof}
    Using \eqref{(21)} we have
    \begin{align}\label{2.43}
      2 \pair{\ricc, \Delta \ricc} &= 2 R_{ik}\Delta R_{ik} = \\ \nonumber &= 2\lambda\abs{\ricc}^2 + 2 \operatorname{Tr}(\ricc^3) + \frac 12 \pair{\nabla\abs{\ricc}^2, X } + R_{ik}S_{ik} - \\ \nonumber &-X_{js}R_{skji}R_{ik} - X_sR_{ik}R_{is, k} - 4 R_{ik}R_{sj}R_{skji} + \\ \nonumber &+X_{jkij}R_{ik}-X_{ikjj}R_{ik}.
    \end{align}
    First we analyze the term $X_{jkij}R_{ik}$. Towards this aim we consider the soliton equation
    \[\tag{\ref{(1)componentsSolEq}}
      R_{ij} + \frac 12 \pa{X_{ij} + X_{ji}} = \lambda \delta_{ij}.
    \]
    Tracing with respect to $i$ and $j$ we obtain
    \[
    S + X_{tt} = m\lambda,
    \]
    so that, taking covariant derivatives,
    \begin{equation}\label{2.44_cov_derTrace_sol_eq}
      S_i = -X_{tti}
    \end{equation}
    and similarly from \eqref{2.44_cov_derTrace_sol_eq}
     \[
      S_{ik} = -X_{ttik}.
    \]
    It follows that
    \begin{equation}\label{2.45}
      R_{ik}S_{ik} = - X_{ttik}R_{ik}.
    \end{equation}
    From the commutation relations \eqref{2.19} and \eqref{2.17covDerivY3rdComm} we get
    \begin{align*}
      X_{jkij} &= X_{jkji} + X_{tk}R_{ti} + X_{jt}R_{tkij} = \\ &=X_{jjki} + X_{si}R_{sk}+X_{s}R_{sk, i} + X_{tk}R_{ti}+X_{jt}R_{tkij}
    \end{align*}
    and therefore, using \eqref{2.45} and soliton equation \eqref{(1)componentsSolEq},
    \begin{align*}
      R_{ik}X_{jkij} &= -S_{ik}R_{ik} + R_{ik}\pa{X_{kt}+X_{tk}}R_{ti} + X_{s}R_{ik}R_{sk, i} + \\ &+X_{jt}R_{ik}R_{tkij}= \\ &= - S_{ik}R_{ik} + 2\lambda\abs{\ricc}^2 -2\operatorname{Tr}(\ricc^3) + X_{s}R_{ik}R_{sk, i} + \\ &+ X_{jt}R_{ik}R_{tkij}.
    \end{align*}
    Substituting this latter into \eqref{2.43} and simplifying we obtain
    \begin{align}\label{2.46}
      2R_{ik}\Delta R_{ik} &= 4 \lambda \abs{\ricc}^2 + \frac 12 \pair{\nabla \abs{\ricc}^2, X} - 2 X_{js}R_{ik}R_{skji}- \\ \nonumber &-4R_{ik}R_{sj}R_{skji} -X_{iktt}R_{ik}.
    \end{align}
    Now we analyze the term $X_{iktt}R_{ik}$. Towards this end we take covariant derivative of the soliton equation \eqref{(1)componentsSolEq}:
    \[
    R_{ij, k} = -\frac 12 \pa{X_{ijk} + X_{jik}}.
    \]
    Tracing with respect to $j$ and $k$ we get
    \[
    R_{ik, k} = -\frac 12 \pa{X_{ikk} + X_{kik}},
    \]
    so that, using \eqref{(16)}, \eqref{2.4'} and \eqref{2.44_cov_derTrace_sol_eq},
    \[
    S_k = -X_{ktt}-X_{tkt}=-X_{ktt}-X_{ttk}-X_{s}R_{sk}=S_k-X_{ktt}-X_{s}R_{sk},
    \]
    that is,
    \[
    X_{itt}=-X_{s}R_{si}.
    \]
    Taking covariant derivative of this latter
    \begin{equation}\label{2.47_cov_der_X_itt}
      X_{ittk}=-X_{sk}R_{si}-X_{s}R_{si, k}.
    \end{equation}
    Next, from \eqref{2.19} and \eqref{2.17covDerivY3rdComm} we obtain
    \begin{align*}
      X_{ittk} &= X_{itkt}+X_{st}R_{sitk}+X_{is}R_{sttk}= \\ &=X_{iktt} + 2X_{st}R_{sitk}-X_{is}R_{sk}+X_{s}R_{sitk, t}.
    \end{align*}
    Hence, using \eqref{2.47_cov_der_X_itt} and \eqref{2.22Ricci23} we deduce
    \begin{align}\label{2.48}
      R_{ik}X_{iktt} &= -X_{s}R_{ik}R_{si, k} - X_{s}R_{ik}R_{sitk, t} -2X_{st}R_{sitk}R_{ik} = \\ \nonumber &= -X_{s}R_{ik}R_{si, k}-2X_{st}R_{sitk}R_{ik}+X_sR_{ik}\pa{R_{ks, i}-R_{ki, s}}= \\ \nonumber &= - \frac 12 \pair{\nabla \abs{\ricc}^2, X} -2X_{st}R_{sitk}R_{ik}.
    \end{align}
    We substitute \eqref{2.48} into \eqref{2.46} to get
    \begin{equation}\label{2.49}
      2R_{ik}\Delta R_{ik} = 4\lambda\abs{\ricc}^2 + \pair{\nabla \abs{\ricc}^2, X} + 4R_{ik}R_{sj}R_{skij}.
    \end{equation}
    Thus, from \eqref{2.40_Delta_abssquared_T}, \eqref{2.49} and \eqref{2.25DeltaScalar} we obtain
    \begin{align}\label{2.50_DeltaabsTsq}
      \Delta \abs{T}^2 &= 2 \abs{\nabla\ricc}^2 - \frac 2m \abs{\nabla S}^2 + 4\lambda \abs{\ricc}^2 + \pair{\nabla \abs{\ricc}^2, X} + \\ \nonumber &+ 4R_{ik}R_{sj}R_{skij}-4 \frac{\lambda}{m}S^2 -\frac 2m S \pair{\nabla S, X} + \frac{4}{m^2}S^3 + \\ \nonumber &+ \frac 4m S \abs{T}^2.
    \end{align}
    An immediate computation shows that
    \[
    \abs{\nabla T}^2 = \abs{\nabla \ricc}^2 -\frac 1m \abs{\nabla S}^2.
    \]
    Using this fact and \eqref{2.39_traceless_Ricci_norm}, after some algebraic manipulation from \eqref{2.50_DeltaabsTsq} we obtain \eqref{2.42}.
  \end{proof}
\end{proposiz}
We recall the decomposition of the curvature tensor $(m \geq 3)$ into its irreducible components:
\begin{align}\label{2.51_comp_curvat_tensor}
  R_{ijks} &= W_{ijks} + \frac{1}{m-2}\pa{R_{ik}\delta_{js}-R_{is}\delta_{jk}+R_{js}\delta_{ik}-R_{jk}\delta_{is}}- \\ \nonumber &-\frac{S}{(m-1)(m-2)}\pa{\delta_{ik}\delta_{js}-\delta_{is}\delta_{jk}},
\end{align}
where $W_{ijks}$ are the component of the Weyl curvature tensor $W$. Note that $\varr$, $m=\dim M \geq 3$, is conformally flat if and only if $W \equiv 0$.

We are now ready to prove
\begin{coroll}\label{CO_2.52_Delta_absTsquared}
  Let $\soliton$ be a Ricci soliton with soliton constant $\lambda$ on $\varr$ and let $S(x)$ be the scalar curvature. Assume $m\geq 3$ and that $\varr$ is conformally flat. Then
  \begin{align}\label{2.53_Delta_absTsquared}
    \frac 12 \Delta \abs{T}^2 &= \frac 12 \pair{\nabla \abs{T}^2, X} + \abs{\nabla T}^2 + 2\pa{\lambda - \frac{m-2}{m(m-1)}S}\abs{T}^2 + \\ \nonumber &+\frac{4}{m-2} \operatorname{Tr}(T^3).
  \end{align}
  \begin{proof}
    A computation shows that
    \begin{equation}\label{2.54}
      \operatorname{Tr}(T^3) = \operatorname{Tr}(\ricc^3) - \frac 3m S \abs{\ricc}^2 + \frac{2}{m^2}S^3.
    \end{equation}
    Thus a simple algebraic manipulation using $W \equiv 0$, \eqref{2.51_comp_curvat_tensor} and \eqref{2.54} gives
    \begin{equation}\label{2.55}
      2R_{ik}R_{sj}R_{skij} = \frac{4}{m-2}\operatorname{Tr}(T^3) - 2 \frac{2m-3}{m(m-1)}S\abs{T}^2 - \frac{2}{m^2}S^3.
    \end{equation}
    Inserting \eqref{2.55} into \eqref{2.42} immediately yields the desired equation \eqref{2.53_Delta_absTsquared}.
      \end{proof}
\end{coroll}
\begin{oss}
  In case $X = \nabla f$, $f \in \cinf$, that is, the soliton is a gradient soliton, we can rewrite \eqref{2.53_Delta_absTsquared} in the form
  \[
  \frac 12 \Delta_f \abs{T}^2 =  \abs{\nabla T}^2 + 2\pa{\lambda - \frac{m-2}{m(m-1)}S}\abs{T}^2 +\frac{4}{m-2} \operatorname{Tr}(T^3),
  \]
  which is formula (2.21) (with $\lambda$ constant) of Corollary 2.7 in \cite{prrims}.
\end{oss}

\end{section}

\begin{section}{Proof of the main results}

\begin{subsection}{Proof of Theorem \ref{TheoremANonextLp} and a further result}

First of all, using Cauchy-Schwarz inequality we have that for any vector field $Y$ on $M$
\begin{equation}\label{3.1CauchySchwarz}
   \frac 14 \abs{\nabla \abs{Y}^2}^2 \leq \abs{Y}^2\abs{\nabla Y}^2.
\end{equation}
We set $u = \abs{X}^2$, we multiply \eqref{2.9BochnerNGSol} by $u$ and use \eqref{3.1CauchySchwarz} to obtain
\begin{equation}\label{3.1'}
\frac 12 u\Delta u + u\ricc\pa{X, X} \geq \frac 14 \abs{\nabla u}^2.
\end{equation}
Next we use assumption \eqref{1.19CondRiccileq1.3} to deduce
\begin{equation}\label{3.2inequDeltau}
  u\Delta u + a(x)u^2 \geq \frac 12 \abs{\nabla u}^2.
\end{equation}
From the work of Fisher-Colbrie and Schoen, \cite{fcs}, we know that assumption \eqref{1.4newFirstEigenLgeq0} implies the existence of $\varphi \in C^2(M)$, $\varphi >0$, solution of
\begin{equation}\label{3.3}
  \Delta\varphi + Ha(x)\varphi=0 \qquad \text{on } M.
\end{equation}
Next, we apply the proof of Theorem 3.1 of \cite{prslogtype} with $b(x) \equiv 0$, $K =0$, $A= -\frac 12$ under assumption \eqref{1.5newIntConditionX} which replaces assumption $(3.6)$ of Theorem 3.1 with $p=2$ and $0 \leq \beta \leq H-1$, to arrive up to the conclusion
\begin{equation}\label{3.4uHvarphi}
  u^H = C \varphi
\end{equation}
for some $C \geq 0$. Since by assumption $X \not \equiv 0$ we conclude that $C>0$ and $u>0$ on $M$. We insert the expression of $\varphi$ in terms of $u$ in \eqref{3.3} and divide by $Hu^{H-2}$ to obtain
\[
u\Delta u + a(x)u^2 = -(H-1)\abs{\nabla u}^2.
\]
Thus, from \eqref{3.2inequDeltau} we deduce that $u$ and hence $\abs{X}^2$ are constant. We then go back to \eqref{3.1'} to obtain, using \eqref{1.19CondRiccileq1.3} and \eqref{2.9BochnerNGSol},
\begin{equation}\label{3.5riccnablaX}
  \frac 12 a(x) \abs{X}^2 \geq \ricc\pa{X, X} = \abs{\nabla X}^2.
\end{equation}
However, from \eqref{3.4uHvarphi} $\varphi$ is a positive constant and \eqref{3.3} implies, since $H \geq 1$, $a(x) \equiv 0$. Thus from \eqref{3.5riccnablaX} $\abs{\nabla X} \equiv 0$ on $M$ and $X$ is a parallel vector field. Thus $X$ is a Killing field and going back to \eqref{1.1soliton_def}
\[
\ricc = \lambda \metric,
\]
that is, $\varr$ is Einstein.
Now, since $X$ is parallel, $X$ is a closed conformal field and the final part of the Theorem follows from $(c)$ of Proposition 2 in \cite{montiel} and from Corollary 9.107 of \cite{Besse}. \hfill \qedsymbol

The next result is a consequence of the previous proof.
\begin{proposiz}\label{PR3.6Lpcond}
  In the assumptions of Theorem \ref{TheoremANonextLp} suppose $\operatorname{vol}(M) = +\infty$. Thus, there are no soliton structures $\soliton$ on $\varr$ with $X \not \equiv 0$ and $X \in L^p(M)$ for some $p>0$.
  \begin{proof}
    We proceed as above up to showing that $\abs{X}^2$ is a positive constant. Thus the result immediately follows since $\operatorname{vol}(M) = +\infty$.
  \end{proof}
\end{proposiz}
\end{subsection}

\begin{subsection}{Proof of Theorem \ref{TheoremA'SobolevNonextLp}}

  Let $\abs{X}^2 = u$. Then $u$ satisfies \eqref{3.2inequDeltau}. Now apply Theorem 9.12 of \cite{prsbook} to contradict \eqref{1.5''condaplus}. \hfill \qedsymbol
\end{subsection}

\begin{subsection}{Proof of Theorem \ref{TH_triv_expander}}
  From equation \eqref{2.9BochnerNGSol} we deduce
  \begin{equation}\label{DeltaAbsXgeqaX}
    \Delta\abs{X}^2=2\abs{\nabla X}^2 -2\ricc(X, X) \geq \frac{2B}{\pa{1+r(x)}^\mu}\abs{X}^2 \quad \text{on }\, M.
  \end{equation}
  Assume now $\abs{X}^2 \not \equiv 0$ and choose $\gamma >0$ such that $\Omega_\gamma = \set{x \in M : \abs{X}^2 >\gamma} \neq \emptyset$. On $\Omega_\gamma$ we have then, using \eqref{DeltaAbsXgeqaX},
  \begin{equation*}
    \pa{1+r(x)}^\mu\Delta\abs{X}^2 \geq 2B\abs{X}^2 > 2B\gamma >0,
  \end{equation*}
  which implies
   \begin{equation*}
   \inf_{\Omega_\gamma} \pa{1+r(x)}^\mu\Delta\abs{X}^2  > 0.
  \end{equation*}
  Applying now Theorem 4.1 in \cite{prsmax} (with $\varphi(t)=t, A=\delta=1$) we obtain a contradiction. \hfill \qedsymbol
\end{subsection}

\begin{subsection}{Proof of Theorem \ref{TeorBScalarEstimates}}

  From \eqref{2.25DeltaScalar} of Proposition \ref{PR2.24_DeltaScalar} we have
  \begin{equation}\label{3.7}
    \frac 12 \Delta S \leq \frac 12 \pair{\nabla S, X} + \lambda S - \frac {S^2}{m}.
  \end{equation}
  Thus, $u=-S$ satisfies the differential inequality
  \[
  \frac 12 \Delta u \geq \frac 12 \pair{\nabla u, X} +\lambda u + \frac{u^2}{m}.
  \]
  Therefore, from \eqref{eq1.6supX} and the above we have
  \begin{equation}\label{3.8}
    \frac 12 \Delta u \geq -\frac 12 \abs{X}^*\abs{\nabla u} +\lambda u + \frac{u^2}{m}.
  \end{equation}
  Now the validity of the Omori-Yau maximum principle on $M$ implies that of Theorem 1.31 on \cite{prsmax} that we apply with the choices $F(t) = t^2$ and
  \[
  \varphi\pa{u, \abs{\nabla u}} = -\frac 12 \abs{X}^*\abs{\nabla u} +\lambda u + \frac{u^2}{m}.
  \]
  Then $u^* = \sup_M u < +\infty$ and
  \begin{equation}\label{3.9ustarleq0}
    \lambda u^* + \frac{(u^*)^2}{m}\leq 0.
  \end{equation}
  But $u^* = -S_*$, so that the claimed bounds on $S_*$ in the statement of Theorem \ref{TeorBScalarEstimates} follow immediately from \eqref{3.9ustarleq0}.

  \emph{Case (i)}. Suppose now $\lambda <0$ and that for some $x_0 \in M$
  \[
  S(x_0) = S_* = m\lambda.
  \]
  In particular $S(x) \geq m\lambda$ on $M$ and the function $w = S - m\lambda$ is non-negative on $M$. From \eqref{3.7} we immediately see that
  \begin{equation}\label{3.10ineqw}
    \Delta w - \pair{X, \nabla w} +2\lambda w \leq \Delta w -\pair{X, \nabla w} + 2\frac Sm w \leq 0.
  \end{equation}
  We let
  \[
  \Omega_0 = \set{x \in M : w(x) =0}.
  \]
  $\Omega_0$ is closed and not empty since $x_0 \in \Omega_0$; let now $y \in \Omega_0$. By the maximum principle applied to \eqref{3.10ineqw} $w \equiv 0$ in a neighborhood of $y$ so that $\Omega_0$ is open. Thus $\Omega_0 = M$ and $S(x) \equiv \lambda m$ on $M$. From equation \eqref{2.25DeltaScalar} we then deduce $\abs{\ricc - \frac Sm \metric} \equiv 0$, that is, $\varr$ is Einstein and from \eqref{1.1soliton_def} $X$ is a Killing  field. Analogously, if $S(x_0) = S_* = 0$ for some $x_0 \in M$, we deduce that $\varr$ is Ricci flat and $X$ is a homothetic vector field.

  \emph{Case (ii)}. Suppose $\lambda =0$ and that, for some $x_0 \in M$,
  \[
  S(x_0) = S_* =0.
  \]
  From \eqref{3.7}
  \[
  \Delta S - \pair{X, \nabla S} \leq -\frac{S^2}{m} \leq 0.
  \]
  Since $S(x) \geq S_* =0$, by the maximum principle we conclude $S(x) \equiv 0$, by \eqref{2.25DeltaScalar} $\varr$ is Ricci-flat and from \eqref{1.1soliton_def} $X$ is a Killing field.

  \emph{Case (iii)}. Finally, suppose $\lambda>0$. Then $S(x) \geq S_* \geq 0$. From \eqref{3.7}
  \[
  \Delta S - \pair{X, \nabla S} -2\lambda S \leq 0.
  \]
  If $S(x_0) = S_* =0$ for some $x_0\in M$, then again by the maximum principle $S(x)\equiv 0$. From \eqref{2.25DeltaScalar}, $\varr$ is Ricci-flat and from \eqref{1.1soliton_def} $\mathcal{L}_X\metric = 2\lambda\metric$ so that $X$ is a homothetic vector field. Suppose now $S(x_0) =S_* =m\lambda$ for some $x_0 \in M$. From \eqref{3.7}
  \[
  \Delta S - \pair{X, \nabla S} \leq 2\frac Sm \pa{\lambda m-S}
  \]
  and since $S(x) \geq S_* = m\lambda >0$,
   \[
  \Delta S - \pair{X, \nabla S} \leq 0 \qquad \text{on } M.
  \]
  By the maximum principle $S(x) \equiv m\lambda$, from \eqref{2.25DeltaScalar} $\varr$ is Einstein and \eqref{1.1soliton_def} implies that $X$ is a Killing field. Furthermore, since $\lambda >0$, $\varr$ is compact by Myers's Theorem. \hfill \qedsymbol

\end{subsection}
\vspace{0.3 cm}

\begin{subsection}{Proof of Corollary \ref{CO_1.4}}

 Under the curvature assumption of Corollary \ref{CO_1.4}, we have the validity of the Omori-Yau maximum principle (see \cite{prsmax}) and the result now follows from Theorem \ref{TeorBScalarEstimates}. \hfill \qedsymbol

\end{subsection}
\vspace{0.3 cm}

\begin{subsection}{Proof of Corollary \ref{CO_1.4'}}

  Let $\varphi : M \ra \erre^n$ be a proper minimal immersion. Then by \cite{prsmax}, Example 1.14, the Omori-Yau maximum principle holds on $\varr$. If $\varphi$ is not totally geodesic then $S_* = \inf_M S(x) < 0$ and we contradict Theorem \ref{TeorBScalarEstimates}. \hfill \qedsymbol

\end{subsection}
\vspace{0.3 cm}

\begin{subsection}{Proof of Theorem \ref{TH_1.5_Okumura}}

  By Okumura's lemma, \cite{okumura},
  \[
  \operatorname{Tr}(T^3) \geq - \frac{m-2}{\sqrt{m(m-1)}}\abs{T}^3.
  \]
  Thus, from \eqref{2.53_Delta_absTsquared} of Corollary \ref{CO_2.52_Delta_absTsquared} we deduce
  \begin{align*}
    \frac 12 \Delta \abs{T}^2 &\geq \frac 12 \pair{\nabla \abs{T}^2, X} + \abs{\nabla T}^2 + 2\pa{\lambda - \frac{m-2}{m(m-1)}S^*}\abs{T}^2- \\ &- \frac{4}{\sqrt{m(m-1)}}\abs{T}^3 \\ &\geq -\frac 12 \abs{\nabla \abs{T}^2}\abs{X}^*+ 2\pa{\lambda - \frac{m-2}{m(m-1)}S^*}\abs{T}^2 - \frac{4}{\sqrt{m(m-1)}}\abs{T}^3.
  \end{align*}
  Setting $\abs{T}^2=u$ we rewrite the above as
  \begin{equation}\label{u_flower}
    \frac 12 \Delta u \geq -\frac 12 \abs{\nabla u}\abs{X}^* +  2\pa{\lambda - \frac{m-2}{m(m-1)}S^*}u - \frac{4}{\sqrt{m(m-1)}}u^{3/2}.
  \end{equation}
  Now if $\abs{T}^*=+\infty$ then \eqref{1.11_Okumura_est} is obviously satisfied. Otherwise $u^* = \sup_M u < +\infty$ and we can apply the Omori-Yau maximum principle to \eqref{u_flower} to obtain
  \[
  \frac{4}{\sqrt{m(m-1)}}u^*\sq{\frac 12 \pa{\lambda\sqrt{m(m-1)}-\frac{m-2}{\sqrt{m(m-1)}}S^*} - \sqrt{u^*}} \leq 0,
  \]
  from which we deduce that either $u^* = 0$, that is $T \equiv 0$ on $M$, or $\abs{T}^*$ satisfies \eqref{1.11_Okumura_est}. In the first case $\varr$ is Einstein, and being conformally flat it is necessarily of constant sectional curvature. \hfill \qedsymbol

\end{subsection}

\end{section}

\begin{section}{A final remark}
  Let $Y$ be a smooth vector field on $M$. Define the associated vector field
  \[
  W_Y = T\pa{Y, \,}^\sharp,
  \]
  where $T$ is the traceless Ricci tensor and $^\sharp$ is the musical isomorphism $^\sharp : T^*M \ra TM$.
  \begin{lemma}\label{LE_4.1_divWY}
    Let $S(x)$ be the scalar curvature of $\varr$. Then
    \begin{equation}\label{4.2_divWY}
      \diver W_Y = \frac 12 \operatorname{Tr}\pa{\mathcal{L}_Y\metric \circ T} + \frac{m-2}{2m}Y(S).
    \end{equation}
    \begin{proof}
      We give the short proof for completeness. With the notations of Section 2,
      \[
      W_Y = Y_iT_{ij}e_j, \qquad Y=Y_ie_i,
      \]
      where $\set{e_i}$ is the o.n. frame dual to $\set{\theta^i}$. Thus
      \[
      \diver W_Y = Y_{ik}T_{ik} + Y_iT_{ik, k}.
      \]
      Using the fact that $T$ is symmetric and \eqref{(16)} we have
      \begin{align*}
        \diver W_Y &= \frac 12 \pa{Y_{ik}+Y_{ki}}T_{ik} + Y_i\pa{R_{ik, k} - \frac{S_i}{m}} = \\ &=\frac 12 \operatorname{Tr}\pa{\mathcal{L}_Y\metric \circ T} + \frac{m-2}{2m}S_iY_i,
      \end{align*}
      that is, \eqref{4.2_divWY}.
         \end{proof}
  \end{lemma}

      Thus,
      \begin{proposiz}\label{PR_4.2_KazdWartype}
        Let $\soliton$ be a soliton structure on $\varr$ and let $S(x)$ be the scalar curvature. Then
        \begin{equation}\label{4.3_diverWX}
          \diver W_X = \frac{m-2}{2m}X(S) - \abs{T}^2.
        \end{equation}
        In particular, if $M$ is compact
        \begin{equation}\label{4.4_KazdWarnSol}
          \frac{m-2}{2m}\int_M X(S) = \int_M \abs{T}^2.
        \end{equation}
        \begin{proof}
          We use soliton equation \eqref{1.1soliton_def}, \eqref{2.38_traceless_Ricci} and \eqref{2.39_traceless_Ricci_norm} into \eqref{4.2_divWY} to immediately obtain \eqref{4.3_diverWX}.
        \end{proof}
     \end{proposiz}

      \begin{oss}
       Equation \eqref{4.4_KazdWarnSol} can be interpreted as a kind of ``Kazdan-Warner condition'' for compact solitons.
      \end{oss}

\end{section}

\bibliographystyle{plain}
\nocite{*}
\clearpage \addcontentsline{toc}{chapter}{References}
\bibliography{BiblioNGSol}

\end{document}